\documentclass[a4paper,11pt]{article}
\usepackage{latexsym, amssymb,amsmath,amsthm}

\def \To{\longrightarrow}
\def \dim{\operatorname{dim}}
\def \gr{\operatorname{gr}}

\def \C{\mathcal{C}}
\def \D{\Delta}
\def \d{\delta}
\def \e{\varepsilon}
\def \M{\mathrm{M}}
\def \S{\mathcal{S}}
\def \sc{\mathrm{Sc}}

\def \g{\mathrm{g}}

\numberwithin{equation}{section}

\newtheorem{theorem}{Theorem}[section]

\newtheorem{proposition}[theorem]{Proposition}
\newtheorem{corollary}[theorem]{Corollary}
\newtheorem{definition}[theorem]{Definition}

\newtheorem{remark}[theorem]{Remark}

\begin{document}

\title{\large \bf QUIVER APPROACHES TO QUASI-HOPF ALGEBRAS}
\author{Hua-Lin Huang \\ \footnotesize School of Mathematics, Shandong
University, Jinan 250100, P. R. China \\ \footnotesize Email:
hualin@sdu.edu.cn}

\date{}
\maketitle

\begin{abstract}
We provide a quiver setting for quasi-Hopf algebras, generalizing
the Hopf quiver theory. As applications we obtain some general
structure theorems, in particular the quasi-Hopf analogue of the
Cartier theorem and the Cartier-Gabriel decomposition theorem.
\end{abstract}

\section{Introduction}

The notion of quasi-Hopf algebras was introduced by Drinfeld
\cite{d} in connection with the Knizhnik-Zamolodchikov system of
equations. It is obtained from that of Hopf algebras by a weakening
of the coassociativity axiom. Quasi-Hopf algebras turn out to be
very useful in various areas of mathematics and physics such as
low-dimensional topology, number theory, integrable systems, and
conformal field theory.

Quivers are oriented graphs consisting of vertices and arrows.
They are widely used in many areas of mathematics and physics. In
particular thanks to their combinatorial behavior, quivers are
very powerful in the investigation of algebraic structures and
representation theory.

We propose to carry out a systematic study of elementary quasi-Hopf
algebras and pointed dual quasi-Hopf algebras \cite{majid1} by
taking advantage of quiver techniques (see e.g. \cite{ars}). The
goal of the present paper is to provide a handy quiver setting. For
a wider setting and the convenience of exposition, we work mainly
with dual quasi-Hopf algebras. A standard dualisation process will
give the corresponding results for quasi-Hopf algebras. To avoid too
many dual's and quasi's we use the term ``Majid algebra" for ``dual
quasi-Hopf algebra", which was proposed by Shnider-Sternberg
\cite{ss1}.

Throughout, we work over a field $k.$ We show in Section 3 that the
path coalgebra $kQ$ of a quiver $Q$ admits a Majid algebra structure
if and only if $Q$ is a Hopf quiver \cite{cr2}, and that any
coradically graded pointed Majid algebra $H$ can be embedded into a
Majid algebra structure on the path coalgebra of some unique Hopf
quiver $Q(H)$ determined completely by $H.$ This generalizes the
quiver setting for Hopf algebras \cite{cr1,cr2, gs, voz} into the
broader class of quasi-Hopf algebras. As applications we obtain some
general structure theorems in Section 4, namely the quasi-Hopf
analogue of the Cartier theorem and the Cartier-Gabriel
decomposition theorem for Hopf algebras (see e.g.
\cite{cartier,mont1,sweedler}). In particular we show that a
cocommutative connected Majid algebra over a field of characteristic
zero is isomorphic to the universal enveloping algebra of a Lie
algebra, which indicates that there is no cocommutative connected
Majid algebra out of the usual Hopf setting.

\section{Majid Algebras and Hopf Quivers}

Introduction to Hopf algebras and Majid algebras can be found in the
books \cite{sweedler, majid2, kassel, ss1}. For basic knowledge of
quivers and their applications to algebras and representation theory
see \cite{ars}.

\subsection{Majid Algebras}

A dual quasi-bialgebra, or Majid bialgebra for short, is a coalgebra
$(H,\D,\e)$ equipped with a compatible quasi-algebra structure.
Namely, there exist two coalgebra homomorphisms $$\M: H \otimes H
\To H, \ a \otimes b \mapsto ab, \quad \mu: k \To H,\ \lambda
\mapsto \lambda 1_H$$ and a convolution-invertible map $\Phi:
H^{\otimes 3} \To k$ called reassociator, such that for all $a,b,c,d
\in H$ the following equalities hold:
\begin{eqnarray}
a_1(b_1c_1)\Phi(a_2,b_2,c_2)=\Phi(a_1,b_1,c_1)(a_2b_2)c_2,\\
1_Ha=a=a1_H, \\
\Phi(a_1,b_1,c_1d_1)\Phi(a_2b_2,c_2,d_2)=\Phi(b_1,c_1,d_1)\Phi(a_1,b_2c_2,d_2)\Phi(a_2,b_3,c_3),\\
\Phi(a,1_H,b)=\e(a)\e(b).
\end{eqnarray}
Here and below we use the Sweedler sigma notation $\D(a)=a_1 \otimes
a_2$ for the coproduct. $H$ is called a Majid algebra if, moreover,
there exist a coalgebra antimorphism $\S: H \To H$ and two
functionals $\alpha,\beta: H \To k$ such that for all $a \in H,$
\begin{eqnarray}
\S(a_1)\alpha(a_2)a_3=\alpha(a)1_H, \quad
a_1\beta(a_2)\S(a_3)=\beta(a)1_H, \\
\Phi(a_1,\S(a_3),
a_5)\beta(a_2)\alpha(a_4)=\Phi^{-1}(\S(a_1),a_3,\S(a_5))
\alpha(a_2)\beta(a_4)=\e(a).
\end{eqnarray}

A Majid algebra $H$ is said to be pointed, if the underlying
coalgebra is pointed. That is, all the simple subcoalgebras of $H$
are one-dimensional. For a given pointed Majid algebra $(H,\D, \e,
\M, \mu, \Phi,\S,\alpha,\beta),$ let $\{H_n\}_{n \ge 0}$ be its
coradical filtration, and $\gr H = H_0 \oplus H_1/H_0 \oplus H_2/H_1
\oplus \cdots$ the corresponding coradically graded coalgebra. It is
routine to verify that $\gr H$ has an induced Majid algebra
structure similar to the Hopf case in \cite{mont1}. The
corresponding graded reassociator $\gr\Phi$ satisfies
$\gr\Phi(\bar{a},\bar{b},\bar{c})=0$ for all
$\bar{a},\bar{b},\bar{c} \in \gr H$ unless they all lie in $H_0.$
Similar condition holds for $\gr\alpha$ and $\gr\beta.$ In
particular, $H_0$ is a sub Majid algebra and turns out to be the
group algebra $kG$ of the group $G=G(H),$ the set of group-like
elements of $H.$

\subsection{Hopf Quivers}

A quiver is a quadruple $Q=(Q_0,Q_1,s,t),$ where $Q_0$ is the set of
vertices, $Q_1$ is the set of arrows, and $s,t:\ Q_1 \longrightarrow
Q_0$ are two maps assigning respectively the source and the target
for each arrow. A path of length $l \ge 1$ in the quiver $Q$ is a
finitely ordered sequence of $l$ arrows $a_l \cdots a_1$ such that
$s(a_{i+1})=t(a_i)$ for $1 \le i \le l-1.$ By convention a vertex is
said to be a trivial path of length $0.$ The path coalgebra $kQ$ is
the $k$-space spanned by the paths of $Q$ with counit and
comultiplication maps defined by $\e(g)=1, \ \D(g)=g \otimes g$ for
each $g \in Q_0,$ and for each nontrivial path $p=a_n \cdots a_1, \
\e(p)=0,$
\begin{equation}
\D(a_n \cdots a_1)=p \otimes s(a_1) + \sum_{i=1}^{n-1}a_n \cdots
a_{i+1} \otimes a_i \cdots a_1 \nonumber + t(a_n) \otimes p \ .
\end{equation}
The length of paths gives a natural gradation to the path coalgebra.
Let $Q_n$ denote the set of paths of length $n$ in $Q,$ then
$kQ=\oplus_{n \ge 0} kQ_n$ and $\D(kQ_n) \subseteq
\oplus_{n=i+j}kQ_i \otimes kQ_j.$ Clearly $kQ$ is pointed with the
set of group-likes $G(kQ)=Q_0,$ and has the following coradical
filtration $$ kQ_0 \subseteq kQ_0 \oplus kQ_1 \subseteq kQ_0 \oplus
kQ_1 \oplus kQ_2 \subseteq \cdots.$$ Hence $kQ$ is coradically
graded. The path coalgebras can be presented as cotensor coalgebras,
so they are cofree in the category of pointed coalgebras and enjoy a
universal mapping property (see e.g. \cite{voz}).

According to \cite{cr2}, a quiver $Q$ is said to be a Hopf quiver
if the corresponding path coalgebra $kQ$ admits a graded Hopf
algebra structure. Hopf quivers can be determined by ramification
data of groups. Let $G$ be a group, $\C$ the set of conjugacy
classes. A ramification datum $R$ of the group $G$ is a formal sum
$\sum_{C \in \C}R_CC$ of conjugacy classes with coefficients in
$\mathbb{N}=\{0,1,2,\cdots\}.$ The corresponding Hopf quiver
$Q=Q(G,R)$ is defined as follows: the set of vertices $Q_0$ is
$G,$ and for each $x \in G$ and $c \in C,$ there are $R_C$ arrows
going from $x$ to $cx.$ For a given Hopf quiver $Q,$ the set of
graded Hopf structures on $kQ$ is in one-to-one correspondence
with the set of $kQ_0$-Hopf bimodule structures on $kQ_1.$

\section{Quiver Setting for Majid Algebras}

A Majid algebra $H$ is a priori a coalgebra. If $H$ is pointed, then
by \cite{cm} there is a unique quiver $Q(H)$ such that $H$ can be
viewed as a ``large" subcoalgebra of the path coalgebra $kQ(H).$
Here by a ``large" subcoalgebra of $kQ(H)$ is meant it contains at
least the $k$-space $kQ(H)_0 \oplus kQ(H)_1$ spanned by the set of
vertices and the set of arrows. The main aim of this section is to
determine what quivers come up as the quivers of pointed Majid
algebras, and conversely to construct Majid structures from these
quivers.

Firstly we consider for what quiver $Q$ the associated path
coalgebra $kQ$ can be endowed with a graded Majid algebra structure.
It turns out that we have nothing new beyond the Hopf quivers of
Cibils and Rosso \cite{cr2}.

\begin{theorem}
Let $Q$ be a quiver. Then the path coalgebra $kQ$ admits a graded
Majid algebra structure if and only if $Q$ is a Hopf quiver.
\end{theorem}

\noindent{\bf Proof:} We assume first that $Q$ is a Hopf quiver.
Then by \cite{cr2}, there exists a graded Hopf algebra structure on
the path coalgebra $kQ.$ Roughly the graded Hopf structure is
constructed as follows. First by definition there exist a group $G$
and a ramification datum $R$ such that $Q=Q(G,R).$ Then view $kQ_0$
as the group Hopf algebra $kG$ and choose a $kG$-Hopf bimodule
structure on the space $kQ_1.$ Finally the graded Hopf algebra is
obtained by extending the bimodule structure to get a compatible
algebra structure with path coalgebra $kQ$ via the universal mapping
property \cite{voz}. Note that Hopf algebras can be viewed as Majid
algebras with trivial reassociator. Therefore, for the Hopf quiver
$Q,$ the associated path coalgebra $kQ$ admits a fortiori a graded
Majid algebra structure.

Conversely, we prove that if $kQ$ admits a graded Majid algebra
structure then $Q$ is a Hopf quiver. Assume that $(kQ, \D, \e, \M,
\mu, \Phi, \S, \alpha, \beta)$ is a graded Majid algebra. First of
all we can restrict the multiplication $\M$ to $Q_0$ and make it a
group. Indeed, since $\M$ is a coalgebra map we have $\D(gh)=\D(\M(g
\otimes h)=gh \otimes gh$ for all $g,h \in Q_0.$ That is, the
multiplication is closed inside $Q_0.$ For all $g \in Q_0,$ we have
$\alpha(g)\beta(g) \ne 0$ by (2.6) and then $\S(g)g=g\S(g)=1$ by
(2.5). Denote $\S(g)$ by $g^{-1}.$ Note that $g^{-1} \in Q_0$ since
$\S$ is a coalgebra antimorphism. In addition, since $\Phi$ is
convolution-invertible we have $\Phi(f,g,h) \ne 0$ for all $f,g,h
\in Q_0,$ and now by (2.1) we have $f(gh)=(fg)h,$ the associativity.
It follows that endowed with the binary operation $\M$ the set of
vertices $Q_0$ becomes a group. For brevity we denote the group
$(Q_0,\M)$ as $G.$

Let $M$ denote the space $kQ_1.$ Note that $M$ is a $kG$-bicomodule
with structure maps $(\d_{_L},\d_{_R})$ defined by
\begin{gather}
\d_{_L}: \ M \To kG \otimes M, \quad a \mapsto t(a) \otimes a, \\
\d_{_R}: \ M \To M \otimes kG, \quad a \mapsto a \otimes s(a),
\end{gather}
for all $a \in Q_1.$ Let $^gM^h = \{ m \in M \ | \ \d_{_L}(m) = g
\otimes m, \ \d_{_R}(m) = m \otimes h \}$ be the $(g,h)$-isotypic
component of $M.$ Then it is in fact the $k$-span of the arrows with
source $h$ and target $g.$ The graded quasi-algebra structure
induces the following linear $kG$-actions on $M:$
\begin{gather}
\rho_{_L}: kG \otimes M \To M, \quad g \otimes m \mapsto g.m, \\
\rho_{_R}: M \otimes kG \To M, \quad m \otimes g \mapsto m.g,
\end{gather}
where $g.m$ and $m.g$ mean the multiplication $\M(g \otimes m)$ and
$\M(m \otimes g)$ respectively in $kQ.$ It is worthy to make a
remark here that in general $\rho_{_L}$ (resp. $\rho_{_R}$) does not
make $M$ a left (resp. right) $kG$-module, since the associativity
holds only up to a non-zero scalar by (2.1). Precisely, for all
$e,f,g,h \in G$ and $m \in \ ^gM^h$ we have
\begin{gather}
e.(f.m)=\frac{\Phi(e,f,g)}{\Phi(e,f,h)}(ef).m,\\
(m.e).f=\frac{\Phi(h,e,f)}{\Phi(g,e,f)}m.(ef),\\
(e.m).f=\frac{\Phi(e,h,f)}{\Phi(e,g,f)}e.(m.f).
\end{gather}

The axioms of graded Majid algebra ensure that $\rho_{_L}$ and
$\rho_{_R}$ are $kG$-bicomodule morphisms. Namely, since $\M$ is a
coalgebra morphism we have for all $f \in G$ and $m \in \ ^gM^h,$
\begin{gather}
\d_{_L}(f.m)=fg \otimes f.m, \quad \d_{_L}(m.f)=gf \otimes m.f, \\
\d_{_R}(f.m)=f.m \otimes fh, \quad \d_{_R}(m.f)=m.f \otimes hf.
\end{gather}
These equalities lead to $f. \ ^gM^h \ \subseteq \ ^{fg}M^{fh},
\quad ^gM^h .f \ \subseteq \ ^{gf}M^{hf}.$ Notice that $f \in G$ is
invertible, it follows that $f^{-1}. \ ^{fg}M^{fh} \ \subseteq \
^gM^h$ and $^{gf}M^{hf}.f^{-1} \ \subseteq \ ^gM^h,$ and that
$f^{-1}.(f.\ ^gM^h) \ = \ ^gM^h$ and $(^gM^h.f).f^{-1} \ = \ ^gM^h$
by (3.5)-(3.6) and (2.2). It is interesting to note that the
composition of $f^{-1}$ and $f$ actions in the preceding two
identities are given by non-zero scalars
$\frac{\Phi(f^{-1},f,g)}{\Phi(f^{-1},f,h)}$ and
$\frac{\Phi(h,f,f^{-1})}{\Phi(g,f,f^{-1})}$ respectively, which are
not necessarily 1 as usual. By combining the above arguments we have
the following identities of vector spaces
\begin{equation}
f. \ ^gM^h \ =  \ ^{fg}M^{fh}, \quad  ^gM^h .f \ = \ ^{gf}M^{hf} \ .
\end{equation}
In particular it follows that for all $x,g,c \in G,$
$$^{g^{-1}cgx}M^x = \ ^{g^{-1}cg}M^1.x = \ (g^{-1}. \ (^cM^1.g)).x \ .$$
It is clear that $$ \dim_k{ ^{g^{-1}cgx}M^x}=\dim_k{^cM^1} \ .$$ Now
recall the geometric meaning of the isotypic spaces. The above
equation implies that, for all $x \in G$ and all $c' \in C,$ where
$C$ is the conjugacy class containing $c,$ there are $\dim_k{^cM^1}$
arrows going from $x$ to $c'x$ in $Q.$ Let $\mathcal{C}$ be the set
of the conjugacy classes of $G.$ For each $C \in \mathcal{C},$ fix
an element $c \in C$ and set $R_C=\dim_k{^cM^1}.$ Take a
ramification datum of $G$ as $R=\sum_{C \in \mathcal{C}} R_C C.$ It
follows from the previous arguments that $Q$ is exactly the Hopf
quiver $Q(G,R).$ We are done. \hfill$\square$

Next we investigate the construction of (non-trivial) Majid algebras
from given Hopf quivers. Naturally we need a quasi-Hopf analogue of
the notion of Hopf bimodules, whose axioms already appear implicitly
in the preceding proof of Theorem 3.1.

\begin{definition}
Assume that $H$ is a Majid algebra with reassociator $\Phi.$ A
linear space $M$ is called an $H$-Majid bimodule, if $M$ is an
$H$-bicomodule with structure maps $(\d_{_L},\d_{_R}),$ and there
are two $H$-bicomodule morphisms
$$\rho_{_L}: H \otimes M \To M, \ h \otimes m \mapsto h.m, \quad
\rho_{_R}: M \otimes H \To M, \ m \otimes h \mapsto m.h$$ such that
for all $g, h \in H, m \in M,$ the following equalities hold:
\begin{gather}
1_H.m=m=m.1_H,\\
g_1.(h_1.m_0)\Phi(g_2,h_2,m_1)=\Phi(g_1,h_1,m^{-1})(g_2h_2).m^0, \\
m_0.(g_1h_1)\Phi(m_1,g_2,h_2)=\Phi(m^{-1},g_1,h_1)(m^0.g_2).h_2, \\
g_1.(m_0.h_1)\Phi(g_2,m_1,h_2)=\Phi(g_1,m^{-1},h_1)(g_2.m^0).h_2,
\end{gather}
where we use the Sweedler notation $$\d_{_L}(m)=m^{-1} \otimes m^0,
\quad \d_{_R}(m)=m_0 \otimes m_1$$ for comodule structure maps.
\end{definition}

We remark that this notion can be rephrased by that of $H$-bimodules
in the monoidal category of $H$-bicomodules, cf. \cite{majid2}. Of
course, if the 3-cocycle $\Phi$ is trivial then the Majid algebra
$H$ is a usual Hopf algebra and $H$-Majid bimodules are the usual
Hopf bimodules.

Now let $Q$ be a Hopf quiver and assume that
$(H,\D,\e,\M,\mu,\Phi,\S,\alpha,\beta)$ is a graded Majid algebra
structure on the path coalgebra $kQ.$ We have proved that $Q_0$ is a
group, denoted by G, and $kQ_0$ is a sub Majid algebra with
reassociator $\Phi$ and quasi-antipode $(\S,\alpha,\beta).$ For
brevity we denote this Majid algebra by $(kG,\Phi).$ Note that the
restriction of $\Phi$ to $G$ is a 3-cocycle and the space $kQ_1$ is
a natural $(kG,\Phi)$-Majid bimodule with structure maps given by
(3.1)-(3.4).

Conversely, let $Q=Q(G,R)$ be a Hopf quiver and $\Phi$ a 3-cocycle
on the group $G.$ Then $kQ_0=kG$ can be understood as a Majid
algebra with reassociator $\Phi.$ If the space $kQ_1$ can be endowed
with a $(kG,\Phi)$-Majid bimodule structure, then we can construct a
graded Majid algebra structure on the path coalgebra $kQ$ via these
data as follows. The process is similar to \cite{cr2}.

Let $\M_0: \ kQ \otimes kQ \longrightarrow kQ_0$ be the composition
of the canonical projection $\pi_0 \otimes \pi_0:\ kQ \otimes kQ
\longrightarrow kQ_0 \otimes kQ_0$ and the multiplication of the
group algebra $kG=kQ_0,$ and $\M_1: kQ \otimes kQ \longrightarrow
kQ_1$ the composition of the canonical projection
\[ \pi_0 \otimes \pi_1 \bigoplus \pi_1 \otimes
\pi_0: \ kQ \otimes kQ \longrightarrow kQ_0 \otimes kQ_1 \bigoplus
kQ_1 \otimes kQ_0 \] and the sum of left and right quasi
$kG$-actions. Then it is clear that $\M_0$ is a coalgebra map and
$\M_1$ is a $kQ_0$-bicomodule map. Let $\M_n=\M_1^{\otimes n} \circ
\D_2^{(n-1)}: \ kQ \otimes kQ \longrightarrow kQ_n,$ where
$\D_2^{(n-1)}$ is the $(n-1)$-iterated comultiplication of the
tensor product coalgebra $kQ \otimes kQ.$ For any pair of paths $p$
and $q$ with $l(p)+l(q)=m,$ it is easy to see that $\M_n(p \otimes
q)=0$ if $m \ne n.$ Therefore by the universal mapping property
$\M=\sum_{n \ge 0}\M_n: \ kQ \otimes kQ \longrightarrow kQ$ is a
well-defined coalgebra map and moreover respects the length
gradation. The quasi-associativity for the map $\M$ follows from the
quasi-associativity (3.12)-(3.14) of the quasi-bimodule. Note that
the reassociator for $kQ$ is obtained by a trivial extension of
$\Phi$ such that $\Phi(x,y,z)=0$ whenever one of $x,y,z$ lies out of
$kQ_0.$ Hence the map $\M$ defines a quasi-algebra structure and we
get a graded Majid bialgebra structure on $kQ.$ We can also
construct a quasi-antipode $(\mathcal{S},\alpha,\beta)$ again via
the universal mapping property. For all $g \in Q_0,$ recall that
$\S(g)=g^{-1}$ and let $\alpha(g)=1$ and
$\beta(g)=1/\Phi(g,g^{-1},g).$ Extend $\alpha$ and $\beta$ to the
function on $kQ$ by letting $\alpha(p)=\beta(p)=0$ for all
nontrivial paths $p$ in the quiver $Q.$ Let $kQ^{cop}$ denote the
coopposite coalgebra of $kQ.$ Set $\S_0:kQ^{cop} \To kQ_0$ to be the
composition of the projection $\pi_0:kQ^{cop} \To kQ_0$ and the map
of taking inversion of the group $G.$ Set $\S_1:kQ^{cop} \To kQ_1$
to be the composition of the projection $\pi_1:kQ^{cop} \To kQ_1$
and the map $a \mapsto
-\frac{\Phi(s(a),s(a),s(a)^{-1})}{\Phi(t(a),s(a),s(a)^{-1})}(t(a)^{-1}.a).s(a)^{-1}$
for all $a \in Q_1.$ Then $\S_0$ is a coalgebra map and $\S_1$ is a
$kQ_0$-bicomodule map, then there is a coalgebra map $\S:kQ^{cop}
\To kQ$ by the universal mapping property. By direct verification we
can show that it is the desired quasi-antipode map.

We summarize the foregoing arguments in the following:

\begin{proposition}
Let $G$ be a group and $(kG,\Phi,\S,\alpha,\beta)$ a Majid
algebra. Let $Q=Q(G,R)$ be the Hopf quiver associated to a
ramification datum $R$ of $G.$ Then the path coalgebra $kQ$ admits
a graded Majid algebra structure with $kQ_0 \cong
(kG,\Phi,\S,\alpha,\beta)$ as a sub Majid algebra if and only if
$kQ_1$ admits a $(kG,\Phi)$-Majid bimodule structure. Moreover,
the set of such graded Majid algebra structures on the path
coalgebra $kQ$ is in one-to-one correspondence with the set of
$(kG,\Phi)$-Majid bimodule structures on $kQ_1.$
\end{proposition}

We remark that a Hopf quiver can be realized by different groups
with ramification data. So in order to get all the graded Majid
algebra structures on the path coalgebra $kQ$ of a given quiver $Q,$
we should consider all the possible realizations of $Q$ as Hopf
quiver and all the possible $kQ_0$-Majid bimodule structures on
$kQ_1.$ Here comes up a natural question of classifying the category
of $(kG,\Phi)$-Majid bimodules for general $G$ and $\Phi.$ It is
well-known from \cite{rosso,s} that the category of Hopf bimodules
over a finite-dimensional Hopf algebra $H$ is equivalent to the
module category of the Drinfeld double $D(H).$ There is an analogue
for the setting of Majid algebras and Majid bimodules. In
particular, when $G$ is a finite group, the category of
$(kG,\Phi)$-Majid bimodules can be described by the module category
of the twisted quantum double $D^{\Phi}(G)$ introduced in
\cite{dpr}. The author is grateful to the referee for suggesting
that the relation between this reference and the present paper
should be investigated.

Finally we consider general pointed Majid algebras. Let $H$ be a
pointed Majid algebra and $\gr H$ its coradically graded version
as mentioned in subsection 2.1. Write the set $G(H)$ of group-like
elements as $G.$ Then $H_0 \cong (kG,\Phi,\S,\alpha,\beta)$ as
Majid algebras for some appropriate $\Phi,\alpha,\beta$ and
$H_1/H_0$ is a $(kG,\Phi)$-Majid bimodule. Let $Q(H)$ be the
quiver of $H,$ then it must be a Hopf quiver since by construction
$kQ(H)_1 \cong H_1/H_0$ which admits a $kQ_0 \cong
(kG,\Phi)$-Majid bimodule structure. The Gabriel type theorem for
pointed Hopf algebras in \cite{voz} can be generalized to the
following for pointed Majid algebras.

\begin{theorem}
Suppose that $H$ is a pointed Majid algebra and $\gr H$ its graded
version induced by the coradical filtration. Then there is a
unique Hopf quiver $Q(H)$ and a graded Majid algebra structure on
the path coalgebra $kQ(H)$ such that $\gr H$ can be embedded into
it as a sub Majid algebra which contains $kQ(H)_0 \oplus kQ(H)_1.$
\end{theorem}

The proof in \cite{voz} can be modified to the quasi setting, so we
omit the detail. This theorem enables us to construct pointed Majid
algebras exhaustively on Hopf quivers. A classification program of
pointed Majid algebras can be carried out in the quiver framework.
The first step is to classify all $(kG,\Phi)$-Majid bimodules for
general group $G$ and 3-cocycle $\Phi.$ This amounts to a
classification of graded Majid algebras on path coalgebras. This is
achieved in a subsequent work \cite{qha2} by generalizing
\cite{cr1,cr2} further to our Majid setting. The second step is to
classify large sub Majid algebras of those on path coalgebras. This
gives a classification of general pointed coradically graded Majid
algebras. The third step is to carry out a suitable deformation
process (see e.g. \cite{ss2}) to get general pointed Majid algebras
from the graded ones. Certainly the classification problem is very
difficult. In this paper we do not intend to go very far in this
program.

We conclude this section with a corollary of Theorems 3.1 and 3.4.

\begin{corollary}
Let $Q$ be an arbitrary quiver. Then the path coalgebra $kQ$ admits
a Majid algebra structure (not necessarily graded) if and only if
$Q$ is a Hopf quiver.
\end{corollary}

\section{Some Structure Theorems}

We apply the quiver setting to investigate general pointed Majid
algebras. In particular, some structure theorems analogous to the
Cartier theorem and the Cartier-Gabriel decomposition theorem are
obtained.

Let $(H,\Phi,\S)$ be a pointed Majid algebra, $G$ its set of
group-likes and $Q$ its quiver. Assume that $Q$ is the Hopf quiver
$Q(G,R)$ associated to some ramification datum $R=\sum_{C \in
\C}R_CC.$ Note that the quiver $Q$ is connected if and only if the
set $\{c \in C \ | \  R_C \ne 0 \}$ generates the group $G.$ In
general, for each $g \in G$ let $Q(g)$ be the connected component
of $Q$ containing $g.$ Denote by $e$ the unit element of $G.$ The
set $N$ of vertices of $Q(e)$ is a normal subgroup of $G,$ since
it is generated by the set $\{c \in C \ | \  R_C \ne 0 \}$ which
is a union of congacy classes. Each connected component $Q(g)$ is
identical to $Q(e)$ as graphs, and its set of vertices is exactly
the coset $Ng.$ The number of connected components of $Q$ is
exactly the index $[G:N].$

The graphical features of the quiver $Q$ imply similar properties
for $H.$ The coalgebra embedding $H \hookrightarrow kQ$ decomposes
$H$ into blocks in traditional terms of algebra (see
\cite{green}), or link-indecomposable components in the sense of
Montgomery (see \cite{mont2}). Let $H_{(g)}$ be the image of $H$
in $kQ(g),$ i.e., the block (or link-indecomposable component) of
$H$ containing $g.$  We call $H_{(e)}$ the principal block.
Obviously $G(H_{(e)})=N$ and the number of blocks is equal to
$[G:N].$ Moreover, we have the following theorem. It is a
quasi-Hopf analogue of the Theorem 3.2 of \cite{mont2}, which can
be viewed as a generalization of the Cartier-Gabriel decomposition
theorem.

\begin{theorem} Keep the notations as above.
\begin{enumerate}
    \item The map $\operatorname{Tr}_{g}: \ H_{(e)} \longrightarrow
    H_{(g)}$ defined by $p \mapsto pg$ is a coalgebra
    isomorphism.
    \item $H_{(g)}H_{(h)} \subseteq H_{(gh)}$ and $\S(H_{(g)})
    \subseteq H_{(g^{-1})}.$ In particular, $H_{(e)}$ is a Majid
    algebra.
    \item Assume further that $H$ is coradically graded. Then there is a Majid
    algebra isomorphism $H \cong H_{(e)} \#_{\sigma}^\Phi kG/N,$ where
    $\sigma:  G/N \times G/N \longrightarrow N$ is a 2-cocycle and $H_{(e)} \#_{\sigma}^\Phi kG/N$
    is a crossed product twisted by $\Phi.$
\end{enumerate}
\end{theorem}

\noindent{\bf Proof:} The claims (1) and (2) are easy. We only prove
the claim (3). Take a set $T$ of distinct coset representatives of
$N$ in $G.$ In particular, for the unit coset $N$ we take $e$ as its
representative. For any $g \in G,$ write $\overline{g} \in T$ as the
the representative of the coset in that g lies. Then there is a
2-cocycle $\sigma:\ G/N \times G/N \To N$ such that
$\bar{u}\bar{v}=\sigma(\bar{u},\bar{v})\overline{uv}$ for any
$\bar{u},\ \bar{v} \in T.$ It follows by (1) that $H =
\bigoplus_{\bar{u} \in T} H_{(e)}\bar{u} \cong H_{(e)} \otimes
kG/N.$ Since $H$ is coradically graded, by Theorem 3.4 it can be
viewed as a sub Majid algebra of $kQ.$ Hence we can choose a basis
for $H$ consisting of paths or linear combinations of paths with the
same source and target. For any basis elements $p, q \in H_{(e)}$
and $\bar{u},\bar{v} \in T,$ define the operation
$$(p \otimes \bar{u})(q \otimes \bar{v})=p(\bar{u} \triangleright
q) \otimes \Theta(p,q,\bar{u},\bar{v})
\sigma(\bar{u},\bar{v})\overline{uv} \ ,$$ where $\bar{u}
\triangleright q=\bar{u}(q\bar{u}^{-1})$ and
$\Theta(p,q,\bar{u},\bar{v}) \in k$ is equal to
$$\frac{\Phi(s(p),\bar{u},s(q)\bar{v})\Phi(s(q),\bar{u}^{-1},\bar{u}\bar{v})
\Phi(\bar{u},t(q)\bar{u}^{-1},\bar{u}\bar{v})\Phi(t(p),t(q),\bar{u}\bar{v})}
{\Phi(t(p),\bar{u},t(q)\bar{v})\Phi(t(q),\bar{u}^{-1},\bar{u}\bar{v})
\Phi(\bar{u},s(q)\bar{u}^{-1},\bar{u}\bar{v})\Phi(s(p),s(q),\bar{u}\bar{v})}
\ .$$ Here we use $s(p)$ and $t(p)$ to denote respectively the
source and the target of a path $p.$ This operation defines a Majid
algebra structure on the tensor coalgebra $H_{(e)} \otimes kG/N,$
the so-called crossed product twisted by $\Phi.$ Now by direct
calculation we can show that the coalgebra isomorphism
$$H \To H_{(e)} \otimes kG/N, \quad p\bar{u} \mapsto p \otimes
\bar{u}$$ preserves the quasi-algebra structure. Therefore we get
the desired Majid algebra isomorphism $H \cong H_{(e)}
\#_{\sigma}^\Phi kG/N.$ \hfill$\square$

\begin{remark}
Thanks to this theorem, the study of pointed Majid algebras can be
reduced to their principal blocks, or equivalently to the connected
case. We make the technical assumption in (3) to guarantee a simpler
definition and exposition for the crossed product $H_{(e)}
\#_{\sigma}^\Phi kG/N$ twisted by $\Phi.$ One can handle more
general situation by adjusting related works on smash products of
quasi-Hopf algebras (see for instance \cite{bpvo}) to Majid
algebras.
\end{remark}

In the following, we consider cocommutative pointed Majid
algebras. We hope to develop a quasi-Hopf analogue of the Cartier
theorem. Hence from now on, the ground field $k$ is assumed to be
of characteristic zero.

Let $H$ be a cocommutative pointed Majid algebra and $Q$ its quiver.
By the decomposition theorem, we may assume that $H$ is
link-indecomposable, or equivalently the quiver $Q$ is connected. In
this situation, $Q$ must be a multi-loop quiver, that is, a quiver
with only one vertex and with arrows starting and ending at it. If
there were at least two vertices in $Q,$ then there is at least one
arrow such that its source is different from its target. Let $a$ be
such an arrow in $Q.$ Since by \cite{cm} $H$ can be regarded as a
large subcoalgebra of $kQ,$ one may assume that $a \in H.$ This
leads to a contradiction with the cocommutativeness:
$$\D(a)=t(a) \otimes a + a \otimes s(a) \ne a \otimes t(a) + s(a)
\otimes a = \D^{op}(a) \ .$$ Now let $\g$ denote the set $\{x \in H
\ | \ \D(x)= x \otimes 1 + 1 \otimes x \}$ of primitive elements.
For any $x,y,z \in \g,$ their multiplication in $H$ is associative
$$x(yz)=(xy)z$$ according to the axioms (2.1) and (2.4). With
bracket defined by the usual commutator, $\g$ becomes a Lie algebra.
Let $U(\g)$ denote the corresponding universal enveloping algebra.
The classical theorem of Cartier asserts that a cocommutative
connected (=pointed and irreducible) Hopf algebra must be isomorphic
to the universal enveloping algebra of the Lie algebra of its
primitive elements. We show that this is also the case for Majid
algebras. The key point is that, as quasi-algebra a cocommutative
connected Majid algebra is generated by primitive elements.

\begin{theorem}
Let $H$ be a cocommutative connected Majid algebra over a field
$k$ of characteristic zero and $\g$ the set of primitives. Then $H
\cong U(\g).$ In particular, $H$ is a usual Hopf algebra and is
generated by primitives.
\end{theorem}

\noindent{\bf Proof:} By the universal property of enveloping
algebra $U(\g),$ the embedding $\g \hookrightarrow H$ can be
extended to a Majid algebra map $\phi: U(\g) \To H.$ The map $\phi$
is injective, since as coalgebra map its restriction to the space of
primitive elements is injective (see e.g. \cite{mont1}).

On the other hand, as coalgebra $H$ can be embedded into the path
coalgebra of a multi-loop quiver. The space spanned by the loops is
in fact isomorphic to $\g.$ We denote the maximal cocommutative sub
coalgebra of the path coalgebra by $\sc(\g).$ Since $H$ is
cocommutative, it can even be viewed as a sub coalgebra of
$\sc(\g).$

Now we get a series of coalgebra embeddings $U(\g) \hookrightarrow H
\hookrightarrow \sc(\g).$ By the Poincar\'{e}-Birkhoff-Witt theorem
and the structural property of $\sc(\g),$ one can show that the
composition of coalgebra maps actually gives rise to an isomorphism
$U(\g) \cong \sc(\g)$ of coalgebras. This leads to the claimed
isomorphism $H \cong U(\g).$ \hfill$\square$

We remark that the proof is a bit sketchy. The same argument was
used in \cite{h} to provide a simple proof for the classical Cartier
theorem. As corollary of the previous results, cocommutative pointed
Majid algebras must be the usual Hopf algebras with possibly
nontrivial 3-cocycles and isomorphic to the smash product of
universal enveloping algebras and group algebras. In particular,
finite-dimensional cocommutative pointed Majid algebras are finite
group algebras with 3-cocycles.

\section{Summary}

We have built up a quiver setting for pointed Majid algebras. For
the case of elementary quasi-Hopf algebras, that is all its simple
modules are one-dimensional, the quiver setting can be provided
dually by using path algebras of quivers instead of path
coalgebras.

At present there is still a lack of abundant examples and general
structure theorems in the quasi-Hopf algebra theory, let alone the
classification. By taking advantage of quiver techniques, bundles of
examples can be constructed easily on concrete Hopf quivers via the
process showed in Section 3. The graphical information naturally
indicates some general structure theorems. Further, a classification
program may be carried out in the quiver framework.

The quiver setting is expected to be very useful in carrying out a
systematic study for the theory of quasi-Hopf algebras. In
particular, we will show in forthcoming works that the quiver
techniques can help to generalize the celebrated theory of pointed
Hopf algebras (see \cite{as} and references therein) to quasi-Hopf
algebras, and to construct and classify some interesting finite
tensor categories, cf. \cite{eo}.

\vskip 0.7cm

\noindent{\large \bf Acknowledgement:} The author thanks the referee
for valuable comments and suggestions which improve the exposition.
The research was supported by the National NSF of China under grant
number 10601052. Part of the work was done in the Chern Institute of
Mathematics (CIM) supported by the Visiting Scholar Program. The
author thanks CIM, in particular Professor Chengming Bai, for
hospitality.

\end{document}